\documentclass[reqno]{amsart} 

\usepackage{amsmath}
\usepackage{amssymb}
\usepackage{hyperref}
\usepackage{amscd}
\usepackage[all]{xy}
\CompileMatrices
\input{xypic}

\theoremstyle{plain}
\newtheorem{thm}{Theorem}[section]
\newtheorem{prop}[thm]{Proposition}

\newtheorem{lemma}[thm]{Lemma}
\newtheorem{cor}[thm]{Corollary}
\newtheorem{cor*}[thm]{Corollary}
\newtheorem*{thm*}{Theorem}
\newtheorem*{prop*}{Proposition}
\theoremstyle{definition}

\theoremstyle{remark}

\newtheorem{remark}{Remark}

\newcommand\vv{\,\vert\,}
\newcommand\sub{\subseteq}
\newcommand\mf[1]{\mathfrak #1}

\newcommand{\rad}{\operatorname{rad}}

\newcommand\dph{\text{depth}\,}
\newcommand\as{\text{Ass}\,}
\newcommand\rk{\text{rk }}
\newcommand\an{\text{ann}\,}

\newcommand\fg{finitely generated }

\newcommand\res[2]{\text{res}^{#1}_{#2}}

\title{Associated primes in group cohomology and a strange result} 
\author{James A. Schafer}

\date{\today}

\address{\begin{flushleft}\quad Department of Mathematics \\
\quad University of Maryland \\
\quad College Park, Maryland 20742
\end{flushleft}}

\email{jas@umd.edu}

\begin{document}
\maketitle

\begin{abstract} In 2010.  Okuyama \cite{O} proved the following.
 
\begin{thm*} Let $G$ be a finite group and $k$ an algebraically closed field of characteristic $p$ dividing the order of $G$. 
\begin{enumerate}\item Let $\frak p$ be an associated prime of $H^*(G,k)$ of $\dim=s$.  Then there exists an elementary abelian $p$-subgroup $E\sub G$ of rank $s$ with $\frak p=\frak p_E$ and depth$\,H^*(C_GE,k)=s$.
\item Let $E$ be an elementary abelian $p$-subgroup of $G$ of rank $s$ and suppose depth$\,H^*(C_GE,k)=s$.  Then $\frak p_E$ is an associated prime of $H^*(G,k)$.
\end{enumerate}\end{thm*}
 
It is shown this result is true for all fields of characteristic $p$. We also show the following result: 
\begin{thm*} If $E$ is an elementary abelian $p$-subgroup of $G$ and $\dph H^*(G,k)=\rk E$ then up to conjugation $E$ contains the maximal elementary abelian subgroup contained in the center of a Sylow $p$-subgroup of $G$. \end{thm*}

\end{abstract}

\section{Okuyama's result}
\begin{thm*} Let $G$ be a finite group and $k$ a field of characteristic $p$ dividing the order of $G$. 
\begin{enumerate}\item Let $\frak p$ be an associated prime of $H^*(G,k)$ of $\dim=s$.  Then there exists an elementary abelian $p$-subgroup $E\sub G$ of rank $s$ with $\frak p=\frak p_E$ and depth$\,H^*(C_GE,k)=s$.
\item Let $E$ be an elementary abelian $p$-subgroup of $G$ of rank $s$ and suppose depth$\,H^*(C_GE,k)=s$.  Then $\frak p_E$ is an associated prime of $H^*(G,k)$.
\end{enumerate}\end{thm*}
  
If one wishes, by using results of Serre and Wilkerson \cite{W} the first part of the theorem in the general case is easy to show.  Their results show that if $\frak p$ is an associated prime of $H^*(G,k)$ then $\frak p=\frak p_E$ for some elementary abelian $p$-subgroup of $G$ with $\rk E=\dim \frak p$.  By Quillen \cite{Q} there exists $\tau_E\in \frak p_E$ and $\tau_E\notin \frak p_{E'}$ for all $E'$ with $\rk E'=\rk E$. If $\frak p_E=\an\,y$ then $\tau_E\,y=0$ and therefore $\res G{E'}\tau_E\,y=0$.
Since $\res G{E'}\tau_E$ is not in the radical of $H^*(E')$ it is a non-zero divisor and therefore a non-zero divisor in $H^*(C_GE')$.  Hence $\res G{C_GE'}\,y=0$. Now a result of Carlson \cite{C} says $\bar y=\res G{C_GF}y\neq 0$ for some elementary abelian $F$ with $\rk F=\rk E$. Hence $\bar y=\res G{C_GE}\,y\neq 0$. Let $J=\an \bar y$. $\res G{C_GE}(\frak p_E)\sub J$ 
$$\rk E=\dim\frak p_E\geq\dim\,(\res G{C_GE})^{-1}J=\dim J \geq\dim Q\geq\dph H^*(C_GE)\geq \rk E$$
where $Q$ is an associated prime of $H^*(C_GE)$. 

Okuyama's proof in both directions depends on results of Green \cite{G} which are only valid for $p$-groups and a result of Benson \cite{B1} whose proof is questionable.  I asked AI [Claude Opus5.5 (Antropic 2026)] to give a proof of Green's results valid for arbitrary groups and find a proof of Benson's result and then write a paper doing so. The resulting paper is available at http://www.terpconnect.umd.edu/$\sim$jas.  It shows that Okuyama's theorem only needs $k$ to be infinite and not algebraically closed.  This is because of Lemma 6.3 in the paper, which definitely uses $k$ infinite in the proof. Replacing that proof with the proof of the lemma below, valid for any field, allows one to conclude Okuyama's theorem is valid for all fields.

\begin{lemma} Let $Q$ be a finite group acting effectively and linearly on a finite dimensional $k$-vector space $U$ and hence on the symmetric algebra $S(U)$. Then there is $f\in S(U)$ such that $\{{}^qf\,\vert\, q\in Q\}$ is a linearly independent set.\end{lemma}
\begin{proof} Identify $S(U)$ with $k[U^*]$.  The effective action of $Q$ induces an effective action of $Q$ on $K=k(U^*)$. Hence $K$ is a Galois extension of $K^Q$ with Galois group $Q$.  By the Normal Basis theorem there exists $\alpha\in K$ with $\{{}^q \alpha\}$ a $K^Q$-basis of $K$ and therefore linearly independent over $k$.  If $\alpha=\beta/\gamma$ with $\beta,\,\gamma\neq 0$ in $S[U]$, let $\rho=\prod{}^q\gamma\in K^Q$.  Then if  
$\bar\alpha=\alpha\,\rho\in S[U]$, $\{{}^q\bar\alpha\}=\{\rho\,{}^q\alpha\}$ is a linearly independent set in $S(U)$.\end{proof}

\section{A strange result}

Let $G$ be a finite group and $k$ a field of characteristic $p$ dividing the order of $G$.  Let $E_0$ be the unique (up to conjugation) maximal elementary abelian $p$-subgroup of $G$ contained in the center of a Sylow $p$-subgroup of $G$.  Let $C_GE$ denote the centralizer of $E$ in $G$.

\begin{thm*} If $E$ is elementary abelian $p$-subgroup of $G$ and $\dph H^*(C_GE)=\rk E$, then $E_0\sub E$ up to conjugation.\end{thm*}
\begin{proof} By \cite{O}, $\frak p_E$ is an associated prime of $H^*(G)$ and by \ref{fine}, $E$ must contain $E_0$.
\end{proof}

Why is this result strange?  What has the depth of $H^*(C_GE,k)$ have to do with the placement of $E$ in the lattice of elementary subgroups of $G$?

\begin{prop} Let $f:A\to B$ be a ring homomorphism. If $\mf q\sub B$ is a primary ideal with $\mf p=\rad\,\mf q$, then $f^{-1}\mf q\sub A$ is a primary ideal of with $\rad\,f^{-1}\mf q=f^{-1}\mf p$.
\end{prop} 
\begin{proof} If $f(uv)\in\mf q$ with $f(v)\notin\mf q$, then $f(u)^s=f(u^s)\in\mf q$. 
$$ y^s\in f^{-1}\mf q\iff f(y)^s=f(y^s)\in\mf q\iff y\in f^{-1}\rad\mf q.$$\end{proof}

\begin{prop} Let $A\sub B$ be Noetherian.  Then $\as A\sub\{\mf p\bigcap A\,\vert\, \mf p\in\as B\}$.
\end{prop}
\begin{proof}
Since $A$ and $B$ are Noetherian, if $(0)=\bigcap\mf q$ is a minimal primary decomposition of $(0)\sub A$ then $\as A=\{\mf p=\rad q\}$ and similarly for $B$.  Therefore if $(0)=\bigcap\mf q$ is a minimal primary decomposition of $(0)\sub B$ with primes $\{\mf p\}$ then $(0)= \bigcap \mf q\cap A$ is a primary decomposition of $(0)\sub A$ with primes $\{\mf p\cap A\}$ by the above proposition but not necessarily minimal. However it can always be made minimal by a) amalgamating primary ideals with then same radical and b) eleminating a primary ideal $\frak q$ with $\bigcap_{\frak q'\neq\frak q}\sub\frak q$ and therefore $\as A\sub\{\mf p\bigcap A \,\vert\, \mf p\in\as B\}$.

\end{proof}

\begin{thm}\cite{Q} $\frak p_E\sub \frak p_{E'}$ iff $E'$ is conjugate to a subgroup of $E$. \end{thm}

\begin{cor}\ Let \,$E$ be elementary $p$-abelian.  Then there exists $\tau_E\in\frak p_E$ such that $\tau_{E}\notin\frak p_{E'}$ for any elementary abelian $E'$ not conjugate to a subgroup of $E$.\end{cor}
\begin{proof} Let $W=\{E'\in\mathcal A(G)\vv E'\,\text{not conjugate to a subgroup of}\,E\}$. Then $\frak p_E\not\subseteq \frak p_{E'}$ for any $E'\in W$.  By Prime Avoidence $\frak p_E\not\subseteq\bigcup_W\frak p_{E'}$.  Therefore there exists $\tau_E\in\frak p_E$ such that $\tau_{E}\notin\frak p_{E'}$ for any $E'\in W$. \end{proof} 

\begin{remark} Since if $g \in G$ then $\frak p^G_E=\frak p^G_{{}^gE}$ we may always assume the elementary abelian $p$-subgroup in $\frak p^G_E$ has $E\sub P$.
\end{remark}

\begin{cor}\label{fine} Let $E_0$ be the maximal $p$-elementary abelian subgroup of the center of $P$ where $P$ is a Sylow $p$-subgroup of $G$.  Then if $\frak p_E\in\as H^*(G,k)$, $E_0\sub E$.\end{cor}
\begin{proof} Since $H^*(G,k)\sub H^*(P,k)$ the above proposition gives $\as(H^*(G,k))\sub\{\frak p^P_F\cap H^*(G,k)\vv \frak p^P_F\in\as H^*(P,k)\}$.  Hence 
$$\frak p^G_E=\frak p^P_F\cap H^*(G)=\frak p^G_F, \text{ with } \frak p^P_F\sub \as\,H^*(P,k).$$
That is we may always assume that up to conjugation $\frak p^P_E\in\as H^*(P,k)$.  

Let $x\in Z(P)$ be an central element of order $p$ of $P$ and let 
$$\frak p^P_E=\ker\{H^*(P,k)\to H^*(E,k)/\text{Rad}\}\in\as H^*(P,k).$$  Suppose $x\notin E$. Then $E'=(x)$ is not conjugate in $P$ to a subgroup of $E$ and so there exists $\tau\in\frak p_E\sub H^*(P,k)$ with $\tau\notin\frak p^P_{E'}$.  Therefore 
$$0\neq\pi\circ\text{Res}^P_{E'}(\tau)\in H^*(E',k)/\text{Rad}. $$  That is $\text{Res}^P_{E'}(\tau)\in H^*(E',k)$ is not nilpotent.  Since $H^*(E',k)=P[\beta\,\bold u]\otimes \Lambda[\bold u]$ this implies $\text{Res}^P_{E'}(\tau)$ is a non-zero-divisor. Since $E'$ is central in $P$, a result of Duflot, \cite{D}, says $\tau\in H^*(P,k)$ is a non zero-divisor and since $\frak p^P_E=\text{ann}_{H^*(P,k)}(y)$ for some $y\in H^*(P,k)$, $\tau\notin\frak p^P_E$. But $\tau\in \frak p_E$ and therefore every central element of $P$ is contained in $E$ and $E_0\sub E$.
\end{proof}

\end{document}